\documentclass[letterpaper, 10 pt, conference]{ieeeconf}  

\IEEEoverridecommandlockouts                              
\usepackage{bbm}
\usepackage{amsfonts}
\usepackage{mathrsfs}
\usepackage{epsfig}
\usepackage{epsf}
\usepackage{epstopdf}
\usepackage{algorithm}
\usepackage{algpseudocode}
\usepackage{graphics} 
\usepackage{epsfig} 
\usepackage{times} 
\usepackage{amsmath} 
\usepackage{amssymb}  
\usepackage{subfigure}
\usepackage{graphicx}
\usepackage{array,color}
\usepackage{booktabs}
\usepackage{xcolor}
\usepackage{amsmath}
\usepackage{amssymb}
\usepackage[normalem]{ulem}
\usepackage{color}
\newcommand\redsout{\bgroup\markoverwith{\textcolor{red}{\rule[0.55ex]{2pt}{0.4pt}}}\ULon} 
\usepackage{xcolor} 

\newtheorem{theorem}{Theorem}
\newtheorem{proposition}{Proposition}

\newtheorem{lemma}{Lemma}

\newtheorem{definition}{Definition}

\usepackage{multicol}
\usepackage{multirow}
\usepackage{booktabs}

\usepackage{verbatim}

\newcommand{\revision}[1]{\textcolor{black}{#1}}
\newcommand{\cut}[1]{\textcolor{black}{#1}}
\newcommand{\keypoint}[1]{\textcolor{black}{#1}} 
\newcommand{\chen}[1]{\textcolor{black}{#1}}

\title{\LARGE \bf
General formation control for multi-agent systems with double-integrator dynamics
}

\author{Chen~Wang, Weiguo~Xia, Jinan~Sun, Ruifeng~Fan and~Guangming~Xie 
\thanks{This work was supported in part by grants from the National Natural Science Foundation of China (NSFC, No. 61503008, 61603071, 91648120, 61633002, 51575005),
	and the China Postdoctoral Science Foundation (No. 2016T90016, 2015M570013).}
\thanks{C. Wang, R. Fan and G. Xie are with the State Key Laboratory of Turbulence and Complex Systems, {Intelligent Biomimetic Design Lab},
       College of Engineering, Peking University, Beijing 100871, China.
       {\tt\small \{wangchen, frf, xiegming\}@pku.edu.cn}}
\thanks{W. Xia is with the School of Control Science and Engineering, Dalian University of Technology, Dalian 116024, China.
       {\tt\small wgxiaseu@dlut.edu.cn}}
\thanks{J. Sun is with the National Engineering Research Center for Software Engineering, Peking University, Beijing 100871, China.
       {\tt\small sjn@pku.edu.cn}}}

\begin{document}
%

\maketitle
\thispagestyle{empty}
\pagestyle{empty}
%
%
\begin{abstract}
%
We study the \emph{general formation} problem for a group of mobile agents in a plane,
in which the agents are required to maintain a distribution pattern,
as well as to rotate around or remain static relative to a \revision{static/}moving target.
The prescribed distribution pattern is a class of general formations that
the distances between neighboring agents or the distances from each agent to the target do not need to be equal.
Each agent is modeled as a double integrator and can merely perceive the relative information of \revision{the target} and its  neighbors. 
\cut{A distributed control law is designed using the limit-cycle based idea to solve the problem.}
%
One merit of the controller is that it can be implemented by each agent in its Frenet-Serret frame so that only local information is utilized without knowing global information.
%
%
\cut{Theoretical} analysis is provided of the equilibrium of the $N$-agent system
and of the convergence of its converging part.
Numerical simulations are given to show the effectiveness and performance of the proposed controller.
\end{abstract}
%

%

%

\par~

%
\section{Introduction}\label{se:introduction}
%

In recent years, control of multi-agent systems has captured increasing attention
due to both its wide practical potential in various applications,
such as exploration \cite{PeMuXi10}, environmental monitoring \cite{DuMa12}, pursuit and evasion \cite{ChHoIs11,DiYaLi12,Re09}, and surveillance \cite{BuDoHo10},
and its theoretical challenges arising from restrictions  in  application implementations.

Formation control is one of the most actively studied topics within the realm of control and coordination of multi-agent systems,
since in such cooperative tasks the robots can benefit from forming clusters or moving in a desired formation with certain geometric shapes \cite{XiCa11,OhPaAh15,BuCoMa09}.
In particular, by forming desired patterns, the robots are able to successfully complete the tasks \cite{OhPaAh15} and even to improve their performance, such as the quality of the collected data, and the robustness of group motion against random environmental disturbances \cite{BuCoMa09}. 
One theoretical challenge of such formation control problems for multi-agent systems arises from the fact that the robots can use only local information to implement their distributed control strategies without centralized coordination.

Intensive research efforts have been devoted to the distributed formation control for multi-agent \revision{systems}
in the systems and control community \cite{OhPaAh15,AnYuFiHe08,ReBeAt07}.
A considerable amount of studies have focused on consensus based formation control
where the formation control problem is converted by a proper transformation to a state consensus problem.
Specifically, the dynamics of the agents are modeled as single-integrators \revision{\cite{SoLiFe16,WaXi17a}}, double-integrators \cite{ShLiTe15}, and unicycles \cite{MaBrFr04, BrSeCa14, YuLi16a};
some constrained conditions are considered including input saturation \cite{SoLiFe16}, agents' locomotion constraints \cite{WaXiCa12b}, finite-time control \cite{WaXiCa12a} , and limited communication \cut{\cite{GaCaYuAnCa13}.}
%
\cut{With the aid of limit-cycle oscillators, the property of collision avoidance has been guaranteed}
when controlling a group of agents to form a circle around a prescribed target \cite{WaXi17}. Using the nonlinear bifurcation dynamics, including limit cycles, \cite{ChSuLiWa17} has proposed swarm control laws to realize some formation configurations of large-scale \revision{swarms.}
%
%
\cut{From} these studies, the potential of limit-cycle oscillators to formation controllers design has been shown, which greatly inspires our work in this paper.

%
The goal of this paper is to design a distributed controller that can guide a group of mobile agents \revision{with double-integrator dynamics} to form any given general formation in a plane.
The general control objective of the problem comprises two specific sub-objectives.
The first is \emph{target circling} that each agent rotates around or remains \revision{fixed} relative to a \revision{static/}moving target as expected, as well as keeping desired distance to the target.
The second is \emph{distribution adjustment} that each agent maintains the desired distance from its neighbors.
It's worth to emphasize that the \textbf{general formations} allow that the distance between neighboring agents are distinguished
and the distances from \revision{the} agents to the target are different.
%
\revision{In addition, the} agents can only sense local information including the relative information of the target and their two neighbors.

To realize the general formation, a limit-cycle based design is delivered in this paper.
We propose to use a controller comprised of two parts to deal with the two sub-objectives of target circling and distribution adjustment.
The key idea is to first design a limit cycle oscillator as the \emph{converging part}, which makes each agent keep a desired distance to the \revision{static/}moving target as well as  rotating counterclockwise/clockwise around or remaining static relative to the target as required.
Then a \emph{layout part} is introduced to the designed limit cycle oscillator to further make the agents maintain desired distance from its two neighbors.
\revision{Subsequently,} an integrated controller is obtained to solve the general formation problem.
\cut{Our proposed controller can be implemented by agents in their Frenet-Serret frame, so that only local information is utilized without knowing global information.}

The rest of the paper is organized as follows.
In Section \ref{se:problem}, we formulate the general formation problem.
{Then we design a distributed  controller and provide some theoretical analysis on its performances in Section \ref{se:mainresults}.}
Simulation results are given in Section \ref{se:simulations}.
Finally, Section \ref{se:conclusions} concludes this paper.
%

\par~

%
\section{Problem formulation}\label{se:problem}
%
%

%

\begin{figure}[t]
\begin{center}
           \subfigure[Initial states]{\label{fig:InitialStates}
           	\includegraphics[width=0.45\linewidth]{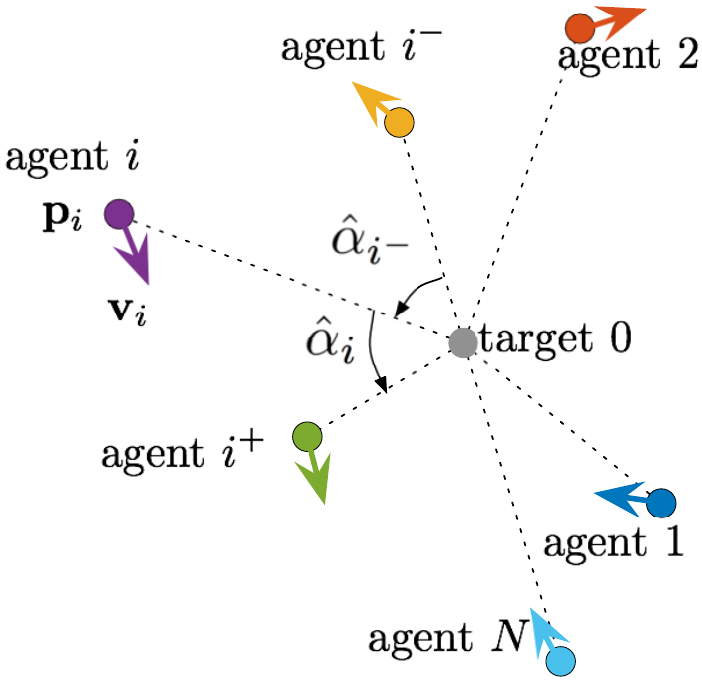}}
           \subfigure[Frenet-Serret frame]{\label{fig:Polar_relative}
           	\includegraphics[width=0.5\linewidth]{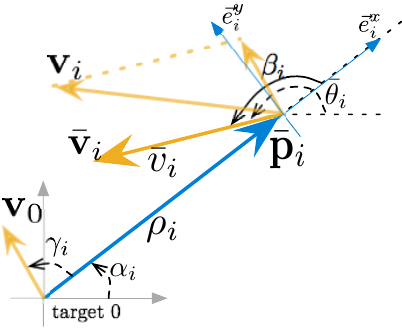}}
   \caption{General formation in a plane. (a) The agents are initially located in a plane.
   		(b) The proposed controller can be implemented in the Frenet-Serret frame of each agent $i$.	
   }\label{fig:Formation}
\end{center}  
\end{figure}

We consider a group of $N$, $N \geq 2$, agents labeled $1$ to $N$ and a \revision{static/}moving target labeled $0$ to be circled around in a plane  (see Fig. \ref{fig:InitialStates}).
%
The $N$ agents' initial positions are not required to be distinguished from each other, whereas no agent occupies the same position as the target.
%
For ease of expression, we label the agents based on their initial positions according to the following three rules:
i) the labels are sorted firstly in ascending order in a counterclockwise manner around the target;
ii) for the agents who lie on the same ray extending from the target, their labels are sorted in ascending order by the distance to the target point;
and iii) for the agents who occupy the same position, their labels are chosen randomly.
Then we consider the case when the agents' neighbor relationships are described
by an undirected ring  graph $\mathbb{G}=(\mathcal{V},\mathcal{E})$,
where $\mathcal{V}=\{1,2,\ldots,N\}$ and $\mathcal{E}=\{(1,2),(2,3),\ldots,(N-1,N),(N,1)\}$.
In such a way, each agent only has two neighbors that are immediately in front of or behind itself.
We denote the set of agent $i$'s two neighbors by $\mathcal{N}_i=\{i^-,i^+\}$ where
\begin{eqnarray}\nonumber
    i^{+} = \begin{cases}
            i+1 \; & \textrm{when} \; i=1, 2, \ldots, N-1\\
            1   \; & \textrm{when} \; i=N
            \end{cases}
\end{eqnarray}\label{eq:neighbor}
and
\begin{eqnarray}\label{eq:topology}
    i^{-} = \begin{cases}
                N   \; & \textrm{when} \; i=1\\
                i-1 \; & \textrm{when} \; i=2,3, \ldots, N.\;\;\;\;\;\;\;\;
            \end{cases}
\end{eqnarray}
%


Let $\mathbf{p}_i  = [x_i, y_i]^T \in \mathbb{R}^2$, $\mathbf{v}_i = [v_i^x, v_i^y]^T \in  \mathbb{R}^2$, and $\mathbf{u}_i = [u_i^x, u_i^y]^T \in \mathbb{R}^2$ denote the position, velocity and control input of agent $i$, respectively.
Each agent $i$ is described by a double-integrator dynamics model
\begin{eqnarray}\label{dynamic:double_i}
	\begin{cases}
		\dot{\mathbf{p}}_i(t) = \mathbf{v}_i(t)\\
		\dot{\mathbf{v}}_i(t) = \mathbf{u}_i(t).
	\end{cases}
\end{eqnarray}
The dynamics of the \revision{static/}moving target are described as follows
\begin{eqnarray}\label{dynamic:double_0}
	\begin{cases}
		\dot{\mathbf{p}}_0(t) = \mathbf{v}_0(t)\\
		\dot{\mathbf{v}}_0(t) = \mathbf{a}_0(t)
	\end{cases}
\end{eqnarray}
where $\mathbf{p}_0  = [x_0, y_0]^T \in \mathbb{R}^2$, $\mathbf{v}_0 = [v_0^x, v_0^y]^T \in  \mathbb{R}^2$, and $\mathbf{a}_0 = [a_0^x, a_0^y]^T\in \mathbb{R}^2$ denote the position, velocity and acceleration of the target, respectively.

\par~



In this paper, the \emph{General Formation} problem in a plane is formalized as to design local controllers for all agents by only using the relative information between the agent and the target and the relative information between the agent and its two neighbors such that all the agents asymptotically form a desired \revision{formation} to keep the \revision{static/}moving target as a reference point.
The general formation is required to rotate clockwise/counterclockwise around the target, or to remain static relative to the target, and to maintain a prescribed distribution pattern \textbf{without} the requirement that all the desired distances between neighboring agents are equal \textbf{nor} the requirement that the desired distances between each agent and the target are equal.


To formulate the problem mathematically, the following variables are introduced.
Let $\bar{\mathbf{p}}_i(t)$ be the relative position between agent $i$ and  target measured by agent $i$ at time $t$,
\begin{equation}\label{eq:bar_p}
	\bar{\mathbf{p}}_i(t) = \mathbf{p}_i(t) - \mathbf{p}_0(t)  \qquad i \in \mathcal{V}
\end{equation}
where $\bar{\mathbf{p}}_i = [\bar{x}_i, \bar{y}_i]^T$.
Denote $\alpha_i(t)$ as the angular of the vector $\bar{\mathbf{p}}_i(t)$ for agent $i$.
The relative velocity between agent $i$ and the target can be derived as
\begin{equation}\label{eq:bar_v}
	\bar{\mathbf{v}}_i(t) = \mathbf{v}_i(t) - \mathbf{v}_0(t)  \qquad i \in \mathcal{V}
\end{equation}
where $\bar{\mathbf{v}}_i = [\bar{v}^x_i, \bar{v}^y_i]^T$.
We further introduce the variables $\hat{\alpha}_i$  as the \emph{angular distance} from agent $i$ to $i^+$,
which is formed by counterclockwise rotating the ray extending from the target to agent $i$ until reaching agent $i^+$.
Similarly, $\hat{\alpha}_{i^{-}}$ is the angular distance from agent $i^-$ to $i$.


Let $d_i$ denote the desired angular spacing from agent $i$ to $i^+$, and $R_i$ denote the desired distance from agent $i$ to the target.
Then the desired distribution pattern of the $N$ agents is determined by the two vectors
\begin{equation}
	d = [d_1,d_2,\ldots,d_N]^T \in \mathbb{R}^N
\end{equation}
and
\begin{equation}
	\mathbf{R} = [R_1,R_2,\ldots,R_N]^T \in \mathbb{R}^N.
\end{equation}
%
Let $\Omega \in \mathbb{R}$ denote each agent's desired angular velocity relative to the target.
For $\Omega > 0 \; (\textrm{resp.} \; \Omega <0)$, the desired formation is required to rotate counterclockwise (resp. clockwise) around the target.
For $\Omega = 0$, the desired formation is required to remain static relative to the target.
{Note that only local information of  $d_i$ and $d_{i^-}$  in vector $d$ is available to each agent $i$.}
We say a prescribed \textbf{\em general formation} $(d_i, R_i, \Omega)$ is {\em admissible} if  $R_i>0$, $d_i>0$
and $\sum^{N}_{i=1}d_i=2\pi$.
%
%



With the above preparation, we are ready to formulate the \emph{General Formation Problem} of interest.

\begin{definition}[General Formation Problem]\label{def:formationG}
	Given an admissible general formation characterized by $\mathbf{R}\in \mathbb{R}^N$ and $d\in \mathbb{R}^N$ in a plane
	with a desired angular velocity $\Omega \in \mathbb{R}$ to a \revision{static/}moving target ${\mathbf{p}}_0$.
	Design distributed control laws $\mathbf{u}_i(t)$, 	$i \in \mathcal{V}$, such that with any initial states
	$[\mathbf{p}_i^T(t_0), \mathbf{v}_i^T(t_0), \mathbf{p}_0^T(t_0), \mathbf{v}_0^T(t_0)] \in\mathbb{R}^8,\forall t_0\geq 0$,
	the solution to system (\ref{dynamic:double_i})
	converges to some equilibrium point $(\mathbf{p}^\ast, \mathbf{v}^\ast)$ satisfying
	\begin{eqnarray}\nonumber
		\| \bar{\mathbf{p}}^*_i\| &=& R_i \\ \label{goal:G1}  
		\dot{\alpha}_i^\ast &=& \Omega  \; \qquad \; i \in \mathcal{V} , \quad  \textrm{\emph{(Target circling)}}\\\nonumber 
		\| \bar{\mathbf{v}}^*_i\| &=& \| R_i \Omega \| 
	\end{eqnarray}
	and
	\begin{eqnarray}\label{goal:G2}
		\hat{\alpha}^* &=& d . \qquad\qquad  \textrm{\emph{(Distribution \revision{adjustment})}}
	\end{eqnarray}
	where ${\cdot}^\ast$ denotes the state at the equilibrium point in this \cut{paper.}
\end{definition}

%

\par~

%
\section{Main results}\label{se:mainresults}
%

In this section, we propose a control law to solve the General Formation \cut{Problem,} 
\cut{and} then give theoretical analysis.

\subsection{Limit-cycle-based control design}



%

\cut{The proposed control law takes the following form:}
\begin{eqnarray}\label{control:part1}
	\mathbf{u}_i = \left[ \begin{array}{cc}  E_i(t) & -\Gamma_i(t) \\ \Gamma_i(t) &  E_i(t) \end{array} \right]  \bar{\mathbf{p}}_i
			+ \left[ \begin{array}{cc} -1 & -1 \\  1 & -1 \end{array} \right] \bar{\mathbf{v}}_i
			+ \mathbf{a}_0
\end{eqnarray}
where
\begin{eqnarray}\label{sys:Polar_para1}{
	\begin{cases}
		E_i(t) = - \mu (\| \bar{\mathbf{p}}_i\| - R_i){\| \bar{\mathbf{p}}_i\|^{\sigma}}  -  \Omega(\Omega - 1)\\
		{\Gamma_i(t)} =  \Omega + {f_i(t)}\\
		{f_i(t)} = \frac{d_{i^{-}}}{d_i+d_{i^{-}}}( \lambda_1 \hat{\alpha}_i + \lambda_2 \dot{\hat{\alpha}}_i) \\
			\quad\quad\quad -\frac{d_i}{d_i+d_{i^{-}}}( \lambda_1 \hat{\alpha}_{i^{-}} +\lambda_2  \dot{\hat{\alpha}}_{i^{-}}),
	\end{cases}}
\end{eqnarray}
and $\lambda_1>0,\lambda_2>0, \mu>0, \sigma \in \mathbb{R} $ are constants.

\cut{Note that the controller is designed in the form of a limit cycle oscillator as the converging part corresponding to the first sub-objective  target circling,
while a layout part $\Gamma_i(t)$ is introduced to deal with the second sub-objective distribution adjustment.}


Let $\rho_i(t) \triangleq \| \bar{\mathbf{p}}_i(t) \|$, $\bar{v}_i(t) \triangleq \| \bar{\mathbf{v}}_i(t) \|$,
and $\beta_i(t)  \triangleq \bar{\theta}_i(t) - \alpha_i(t)$, where $\bar{\theta}_i$ is the angular of the vector ${\bar{v}}_i$.
Then the system (\ref{dynamic:double_i}) under control laws (\ref{control:part1}) can be represented in the polar coordinates
\begin{eqnarray}\nonumber
		\bar{\mathbf{p}}_i(t) &=& \rho_i(t) \left[ \begin{array}{c}  \cos{\alpha_i(t)} \\  \sin{\alpha_i(t)} \end{array} \right] \\\nonumber
		\bar{\mathbf{v}}_i(t) &=& \bar{v}_i(t) \left[ \begin{array}{c}  \cos{\bar{\theta}_i(t)} \\  \sin{\bar{\theta}_i(t)} \end{array} \right]
\end{eqnarray}
as
\begin{eqnarray}\label{sys:Polar}
	\begin{cases}
		\dot{\rho}_i = \bar{v}_i \cos{\beta_i}\\
		\dot{\bar{v}}_i =  { \rho_i (   E_i \cos{\beta_i} + \Gamma_i  \sin{\beta_i}) -  \bar{v}_i } \\
		\dot{\beta}_i = { 1 - \frac{\rho_i}{\bar{v}_i} (   E_i \sin{\beta_i} - \Gamma_i  \cos{\beta_i} ) } - \frac{\bar{v}_i}{\rho_i}\sin{\beta_i}
 	\end{cases}
\end{eqnarray}
and
\begin{eqnarray}\label{sys:Polar_alpha}
	 \dot{\alpha}_i = \frac{\bar{v}_i}{\rho_i}  \sin{\beta_i}
\end{eqnarray}
where $E_i$ and $\Gamma_i$ are given by (\ref{sys:Polar_para1}).

Now we have the overall closed-loop system in the polar coordinates with states ${\rho}_i(t)$, ${\bar{v}}_i(t)$ and ${\beta}_i(t)$ described by equations (\ref{sys:Polar}).
It is worth to emphasize \revision{that the} variables $\alpha_i(t)$ here can be treated as additional states, which are only used for analysis purposes and are not known to the agents (see Fig. \ref{fig:Polar_relative}).


Furthermore, for each agent $i$, we construct a moving frame, the Frenet-Serret frame, that is fixed on the agent with its origin at the representing point and $x$-axis coincident with the orientation of the vector $\bar{\mathbf{p}}_i(t)$.
The agent $i$'s Frenet-Serret frame is shown by $(\vec{e}_i^x, \vec{e}_i^y)$ \cut{in Fig. \ref{fig:Polar_relative}.}
%
\cut{One} can easily check that our proposed control laws (\ref{control:part1}) can be successfully implemented by agents in their Frenet-Serret frame
without knowing the information of global coordinates.

\subsection{Analysis of Equilibrium}



Now, we analyze the equilibria of the $N$-agent system (\ref{dynamic:double_i}) under the control law (\ref{control:part1}).
For this purpose, we consider \revision{both} the closed-loop system (\ref{sys:Polar}) in the polar coordinates and the dynamics of additional states ${\alpha}_i $ described by (\ref{sys:Polar_alpha}).
Then the equilibrium points can be calculated by solving
\begin{eqnarray}\label{cal:equilibrium_pre}
	\begin{cases}
		\dot{\rho}_i = \bar{v}_i \cos{\beta_i} = 0\\
		\dot{\bar{v}}_i = { \rho_i (   E_i \cos{\beta_i} + \Gamma_i  \sin{\beta_i}) -  \bar{v}_i } = 0\\
		\dot{\beta}_i = { 1 - \frac{\rho_i}{\bar{v}_i} (   E_i \sin{\beta_i} - \Gamma_i  \cos{\beta_i} ) } - \frac{\bar{v}_i}{\rho_i}\sin{\beta_i}  =0.
 	\end{cases}
\end{eqnarray}
%


It is known  from the definition of \revision{the} angular distance $\hat{\alpha}_i$ that
\begin{equation}\label{cal:dot_hat_alpha}
	\dot{\hat{\alpha}}_i(t) = \dot{\alpha}_{i^{+}}(t) - \dot{\alpha}_i(t) \qquad i \in \mathcal{V}.
\end{equation}
Together with (\ref{sys:Polar_alpha}), one arrives at a subsystem with states $\hat{\alpha}_i$
%
%
\revision{\begin{eqnarray}\label{dynamic:hat_alpha}
		\ddot{\hat{\alpha}}_i  &=&f_{i^+} -  f_i - \dot{\hat{\alpha}}_i \\\nonumber
					&+& \dot{\alpha}_{i^+}  \cot{\beta_{i^+}}(1-2\dot{\alpha}_{i^+}) \\\nonumber
					&-& \dot{\alpha}_i  \cot{\beta_i}(1-2\dot{\alpha}_i).
\end{eqnarray}}
%


We first analyze the states $\hat{\alpha}, \dot{\hat{\alpha}}$ at the equilibrium point of system (\ref{sys:Polar}). \par~
%
\keypoint{\begin{proposition}\label{pro:hat_alpha_subsystem}
	Any equilibrium point of \revision{the} $N$-agent system (\ref{sys:Polar}) is also an equilibrium of the following system
	%
	%
	\revision{\begin{eqnarray}\label{dynamic:hat_alpha_0}
			\ddot{\hat{\alpha}}_i  &=f_{i^+} -  f_i - \dot{\hat{\alpha}}_i.
	\end{eqnarray}}
\end{proposition}}
%
\begin{proof}
	At any equilibrium point of system (\ref{sys:Polar}),
	one has $\bar{v}^\ast_i = 0$ or $\cos{\beta^\ast_i}=0$ since $\dot{\rho}_i = 0$.
	When $\cos{\beta^\ast_i}=0$, we have $\cot{\beta^\ast_i}=0$.
	When $\bar{v}^\ast_i = 0$ and $\rho^\ast_i \neq 0$,
	we have $\dot{\alpha}^\ast_i = \frac{\bar{v}^\ast_i}{\rho^\ast_i} \sin{\beta^\ast_i}= 0$.
	When $\bar{v}^\ast_i = 0$ and $\rho^\ast_i = 0$,
	it follows that ${\alpha}^\ast_i=0, \bar{\theta}^\ast_i=0$ from their definitions,
	and  $\beta^\ast_i = 0$, and thus $\dot{\alpha}^\ast_i =0$.
	Now, one can conclude that
	$\dot{\alpha}^\ast_i \cot{\beta^\ast_i}  = 0, \; \forall i \in \mathcal{V},$
	always hold at any equilibrium point of system (\ref{sys:Polar}).
	It implies that $\dot{\alpha}^\ast_{i^+}  \cot{\beta^\ast_{i^+}}(1-2\dot{\alpha}^\ast_{i^+}) - \dot{\alpha}^\ast_i  \cot{\beta^\ast_i}(1-2\dot{\alpha}^\ast_i) = 0, \; \forall i \in \mathcal{V}$.
	This completes the proof.
\end{proof}
%

We further rewrite the system (\ref{dynamic:hat_alpha_0}) into a compact form
\begin{equation}\label{dynamic:hat_alpha_0matrix}	
	\left[ \begin{array}{l} \dot{\hat{\alpha}} \\ \ddot{\hat{\alpha}}  \end{array}\right]
	= \Phi(d)
	 \left[\begin{array}{l} {\hat{\alpha}} \\ \dot{\hat{\alpha}}   \end{array}\right]	
\end{equation}
where
\begin{equation}\nonumber	
	\Phi(d)
	= \left[\begin{array}{cc}
		 \mathbf{0}_{N\times N} &  I_N\\
		-\lambda_1 L(d)   &  -\lambda_2 L(d) - I_N
	 \end{array}\right]
\end{equation}
and $ L(d)$ is given by (\ref{eq:Laplacian}).

\begin{figure*}
\begin{equation}\label{eq:Laplacian}
   L(d) = \scriptsize 
   \left[\begin{array}{cccccc}
               \frac{d_2}{d_2+d_1} + \frac{d_N}{d_1+d_N} & - \frac{d_1}{d_2+d_1} & 0 & \ldots & 0 & -\frac{d_1}{d_1+d_N}\\
               - \frac{d_2}{d_2+d_1} & \frac{d_3}{d_3+d_2} + \frac{d_1}{d_2+d_1} & -\frac{d_2}{d_3+d_2} &  \ldots &0 & 0\\
               0 & -\frac{d_3}{d_3+d_2} & \frac{d_4}{d_4+d_3} + \frac{d_2}{d_3+d_2} &\cdots & 0& 0\\
               \vdots & \vdots& \vdots & \vdots &\vdots & \vdots\\
               0 & 0 & 0& \ldots & \frac{d_{N}}{d_N+d_{N-1}} + \frac{d_{N-2}}{d_{N-1}+d_{N-2}} & -\frac{d_{N-1}}{d_N+d_{N-1}}\\
               - \frac{d_N}{d_1+d_N} & 0 & 0 & \ldots  & -\frac{d_{N}}{d_N+d_{N-1}} & \frac{d_1}{d_1+d_N} + \frac{d_{N-1}}{d_N+d_{N-1}}
         \end{array}\right].
\end{equation} 
\end{figure*}


Note that system (\ref{dynamic:hat_alpha_0matrix}), which merely contains variables $\hat{\alpha}, \dot{\hat{\alpha}}$,
is helpful when calculating the equilibrium of the $N$-agent system, especially   {the layout part}  $\Gamma^\ast_i(t)$.
Next we give some useful results about system (\ref{dynamic:hat_alpha_0matrix}) to facilitate the  discussion on the equilibrium of the $N$-agent system (\ref{dynamic:double_i}) under the control law (\ref{control:part1}).

Let $D=\textrm{diag}\{d_1,d_2,\ldots,d_N\}$. Then $D^{-1}L(d)D=L^{T}(d)$.
For analysis purposes, we introduce a pair of  variables $(\delta,\xi)$
\begin{eqnarray}\nonumber
    \begin{cases}
        \delta=D^{-1}\hat{\alpha} \\
        \xi=D^{-1}\dot{\hat{\alpha}}.
    \end{cases}
\end{eqnarray}
Then we  have
\begin{equation}\label{dynamic:hat_alpha_0_tilde}	
	\left[ \begin{array}{l} \dot{\delta} \\ \dot{\xi}  \end{array}\right]
	= \widetilde{\Phi}(d)
	 \left[\begin{array}{l} \delta \\ \xi   \end{array}\right]	
\end{equation}
where
\begin{equation}\nonumber	
	\widetilde{\Phi}(d)
	=\left[\begin{array}{cc}
		 \mathbf{0}_{N\times N} &  I_N\\
		-\lambda_1 L^T(d)   &  -\lambda_2 L^T(d) - I_N
	 \end{array}\right].
\end{equation}
Suppose $\eta_i \;(i=1,2,\ldots,N)$ and $\zeta_{ij}\;(i=1,2,\ldots,N,j=1,2)$
are eigenvalues of $L(d)$ and $\widetilde{\Phi}(d)$, respectively.

%
\begin{lemma}[\keypoint{Lemma 5 of \cite{WaXiCa12a}}]\label{le:eigenvalueL}
   It holds that\\
   \indent   i) $L(d)$ is diagonalizable and $\eta_i\in[0,2]$,     $i=1,2,\ldots,N$;

   ii) $0$ is a single eigenvalue;

   iii) When $N$ is even, $2$ is an eigenvalue, while when $N$ is odd, $2$ is not.
	
\end{lemma}
%

In view of Lemma \ref{le:eigenvalueL}, without loss of generality, we now assume  $\eta_1=0< \eta_2\leq \cdots\leq \eta_N $.
Then we analyze the eigenvalues of  $\widetilde{\Phi}(d)$.
%
\keypoint{\begin{lemma}\label{le:eigenvaluePhi}
    Matrix $\widetilde{\Phi}(d)$ has exactly a zero eigenvalue of algebraic multiplicity $1$
    and all the other eigenvalues have negative real parts.
\end{lemma}}
%
\begin{proof}
Let $\zeta$ be an eigenvalue of the matrix $\widetilde{\Phi}(d)$.
Then, one has $\det(\zeta I_{2N}-\widetilde{\Phi})=0$.
Note that
\begin{eqnarray}\nonumber
    \det(\zeta I_{2N}- \widetilde{\Phi}) = \prod_{i=1}^N [\zeta^2 + (\lambda_2 \eta_i +1) \zeta + \lambda_1 \eta_i] =0.
\end{eqnarray}
Hence,
\begin{eqnarray}\label{eq:eigenvaluePhi}
    \zeta_{i\pm} = \frac{- (\lambda_2 \eta_i +1) \pm \sqrt{ (\lambda_2 \eta_i +1)^2 - 4  \lambda_1 \eta_i}}{2},
    			  \; i \in \mathcal{V}
\end{eqnarray}
From (\ref{eq:eigenvaluePhi}) and Lemma \ref{le:eigenvalueL},
it is easy to see that $\widetilde{\Phi}(d)$ has a zero eigenvalue of algebraic multiplicity $1$
and all the other eigenvalues have negative real parts.
\end{proof}
%
%
\keypoint{\begin{lemma}\label{le:Lemma41RenModified}
    System (\ref{dynamic:hat_alpha_0_tilde}) achieves consensus asymptotically and
     $\delta(t)\rightarrow\mathbf{1}p^T\delta(0)+ \mathbf{1}p^T\xi(0)$, and $\xi(t)\rightarrow \mathbf{0}$,
    as $t$ goes to infinity, where $p\in\mathbb{R}^N$ is the non-negative left eigenvector of $L^T(d)$
    associated with the eigenvalue $0$ and $p^T \mathbf{1}=1$.
\end{lemma}}
%
\begin{proof} 
In view of Lemma \ref{le:eigenvaluePhi},
one can check that eigenvalue zero of $\widetilde{\Phi}(d)$ has geometric multiplicity equal to one.
Note that $\widetilde{\Phi}(d)$ can be written in Jordan canonical form as
\begin{small}
\begin{eqnarray}\nonumber	
	&& \widetilde{\Phi}(d) = PJP^{-1} \\\nonumber
				&=&[\varpi_1, \ldots, \varpi_{2N}]
					\left[\begin{array}{cc} 0 &  \mathbf{0}_{1\times (2N-1)} \\	\mathbf{0}_{(2N-1) \times 1 }   &  J^{\prime} \end{array}\right]
					\left[\begin{array}{c}  \iota_1^T \\ \vdots \\ \iota_{2N}^T\end{array}\right]
\end{eqnarray}
\end{small}
where $\varpi_j$ can be chosen to be the right eigenvectors or generalized eigenvectors of $\widetilde{\Phi}(d)$,
$\iota_j$ can be chosen to be the left eigenvectors or generalized eigenvectors of $\widetilde{\Phi}(d)$,
and $J^{\prime}$ is the Jordan upper diagonal block matrix corresponding to non-zero eigenvalues.

Without loss of generality,  choose $\varpi_1 = [\mathbf{1}_{N}^T, \mathbf{0}_{N}^T]^T$,
where $\mathbf{1}_{N}^T$ and $\mathbf{0}_{N}^T$ are $N$-dimensional all-one and all-zero vectors, respectively. It can be verified that $\varpi_1$ is a right eigenvector of $\widetilde{\Phi}(d)$ associated with the eigenvalue $0$.
%
Let $p$ be the non-negative vector such that $p^T L^T = 0$ and $p^T\mathbf{1} = 1$.
It can be verified that $\iota_1 = [p^T, p^T]^T$ is a left eigenvector of $\widetilde{\Phi}(d)$ associated with  eigenvalue $0$,
where $\iota_1^T \varpi_1 = 1$.

Noting that all eigenvalues of $\widetilde{\Phi}(d)$ except a simple zero eigenvalue have negative real parts, we see that
\begin{eqnarray}\nonumber	
	\mbox{e}^{\widetilde{\Phi}t} &=& P\mbox{e}^{Jt}P^{-1} \\\nonumber
				&=& P
				\begin{small}\left[\begin{array}{cc} 1 &  \mathbf{0}_{1\times (2N-1)} \\	\mathbf{0}_{(2N-1) \times 1 }   &  \mbox{e}^{J^{\prime}t} \end{array}\right]
					P^{-1}\end{small}
\end{eqnarray}
which converges to $\left[\begin{array}{cc} \mathbf{1}p^T &  \mathbf{1}p^T \\ \mathbf{0}_{N\times N}   & \mathbf{0}_{N\times N} \end{array}\right]$ as  $t\rightarrow\infty$.
Noting that
\begin{equation}\nonumber	
	\left[ \begin{array}{l} {\delta}(t) \\ {\xi}(t)  \end{array}\right]
	\rightarrow \left[\begin{array}{cc} \mathbf{1}p^T &  \mathbf{1}p^T \\ \mathbf{0}_{N\times N}   & \mathbf{0}_{N\times N} \end{array}\right]
	 \left[\begin{array}{l} \delta(0) \\ \xi(0)   \end{array}\right]	
\end{equation}
we see that $\delta(t)\rightarrow\mathbf{1}p^T\delta(0)+ \mathbf{1}p^T\xi(0)$, and $\xi(t)\rightarrow \mathbf{0}$ as $t\rightarrow\infty$.
As a result, we know that $|\delta_i(t)-\delta_j(t)| \rightarrow 0$ and $|\xi_i(t)-\xi_j(t)| \rightarrow 0$ as $t \rightarrow \infty$.
That is, system (\ref{dynamic:hat_alpha_0_tilde}) achieves consensus asymptotically.
\end{proof}
%
%
\keypoint{\begin{lemma}\label{th:goal_angleG_0matrix}
	System (\ref{dynamic:hat_alpha_0matrix}) achieves consensus asymptotically.
	Specifically, ${\hat{\alpha}}(t)\rightarrow d$ and $\dot{\hat{\alpha}}(t)\rightarrow \mathbf{0}_N$ as $t \rightarrow \infty$.\end{lemma}}
%
\begin{proof}
From Lemma \ref{le:eigenvaluePhi} and Lemma \ref{le:Lemma41RenModified},
one can see that system (\ref{dynamic:hat_alpha_0_tilde}) achieves consensus asymptotically,
 which further implies that system (\ref{dynamic:hat_alpha_0matrix}) achieves consensus asymptotically.
Moreover, one can check that there exists $p=\frac{1}{2\pi}d$ such that $p^T L^T=0$ and $p^T \mathbf{1}=1$.
Since  $\sum_{i=1}^N {\hat{\alpha}}_i = 2\pi$ and $\sum_{i=1}^N \dot{\hat{\alpha}}_i = 0$ hold all times,
we have
\begin{small}
\begin{eqnarray}\nonumber
    \hat{\alpha} = D\delta & \rightarrow & D\mathbf{1}p^T\delta(0)+ D\mathbf{1}p^T\xi(0)\\\nonumber
    					&=& D\mathbf{1}\frac{1}{2\pi}d^TD^{-1}\hat{\alpha}(0) + D\mathbf{1}\frac{1}{2\pi}d^TD^{-1}\dot{\hat{\alpha}}(0) = d\\\nonumber
    \dot{\hat{\alpha}} = D\xi &\rightarrow& D \mathbf{0} = \mathbf{0}_N.
\end{eqnarray}
\end{small}
Thus, we have $\hat{\alpha} \rightarrow d$ and $\dot{\hat{\alpha}} \rightarrow\mathbf{0}_N$ for large $t$.
\end{proof}
%


With the above preparation, we are ready to calculate the equilibria of the $N$-agent system (\ref{dynamic:double_i}) under the control law (\ref{control:part1}) (i.e., the closed-loop system in the polar coordinates  (\ref{sys:Polar})) by solving (\ref{cal:equilibrium_pre}).
All the equilibrium points can be classified into the following three cases:
\begin{itemize}
	\item[$\bullet$] Case I: $(\rho_i^\ast)^2 + (\bar{v}_i^\ast)^2 \neq 0, \; \forall i \in \mathcal{V}$;
	\item[$\bullet$] Case II: $(\rho_i^\ast)^2 + (\bar{v}_i^\ast)^2 = 0, \; \forall i \in \mathcal{V}$;
	\item[$\bullet$] Case III: $(\rho_i^\ast)^2 + (\bar{v}_i^\ast)^2 \neq 0, \; \exists i \in \mathcal{V}_1$ and $(\rho_j^\ast)^2 + (\bar{v}_j^\ast)^2 = 0, \; \exists j \in \mathcal{V}_2$, where $\mathcal V_1 \bigcup \mathcal V_2= \mathcal{V}$, $\mathcal V_1 \bigcap \mathcal V_2=\emptyset$.
\end{itemize}
%


%
\keypoint{\begin{proposition}[Equilibrium Case I]\label{pro:equilibrium_caseI}
	The equilibrium of the $N$-agent system (\ref{sys:Polar}) is
	(\ref{cal:equilibrium_Ia}) when  $\Omega \neq 0$ and is (\ref{cal:equilibrium_Ib}) when $\Omega = 0$,
	if it satisfies $(\rho_i^\ast)^2 + (\bar{v}_i^\ast)^2 \neq 0, \; \forall i \in \mathcal{V}$.
	\begin{small}
\begin{eqnarray}\label{cal:equilibrium_Ia}
	(\textrm{for} \; \Omega \neq 0) \; { 
		\begin{cases}
			\beta_i^\ast =  \begin{cases}
				\frac{\pi}{2} + 2 k \pi (k \in \mathbb{Z} )    & \textrm{when} \; \Omega >0 \\
				-\frac{\pi}{2} + 2 k \pi (k \in \mathbb{Z} )  & \textrm{when} \; \Omega <0
						\end{cases} \\
			\bar{v}_i^\ast =   \|  \Omega R_i  \|  >0 \\
			\rho_i^\ast = R_i >0  \\
			\dot{\alpha}_i^\ast = \Omega \neq 0	\\
			\hat{\alpha}^\ast = d				
		\end{cases} }
\end{eqnarray}
\end{small}
\begin{small}
\begin{eqnarray}\label{cal:equilibrium_Ib}
	 (\textrm{for} \; \Omega = 0) \;  { 
		\begin{cases}
			\bar{v}_i^\ast =  0 \\
			\rho_i^\ast = R_i >0  \\
			\dot{\alpha}_i^\ast = \Omega =  0	\\
			{\hat{\alpha}}^\ast = d.		
		\end{cases} }
\end{eqnarray}
\end{small}
\end{proposition}}
%
\begin{proof}
In this case, we need to consider  three subcases.

\chen{Subcase I-a: $\bar{v}_i ^\ast \neq 0, \; \forall i \in \mathcal{V}$.}
From (\ref{cal:equilibrium_pre}), one can have $\cos{\beta^\ast_i} = 0$ due to $\bar{v}_i ^\ast \neq 0$, thus $\sin{\beta^\ast_i} = \pm 1$,
and $\bar{v}^\ast_i =  \Gamma_i^\ast \rho^\ast_i \sin{\beta^\ast_i}$ holds.
Together with the definition of $\bar{v}_i \geq 0$ and $\rho_i \geq 0$,
one can check $\bar{v}^\ast_i > 0$ and thus $\rho_i^\ast >0$ and $\Gamma_i^\ast \sin{\beta^\ast_i}>0, \Gamma_i^\ast \neq 0$.
It follows that $\bar{v}_i^\ast =   \|  \Gamma_i^\ast \rho_i^\ast  \| >0$,
and $\sin{\beta_i} = \frac{\bar{v}_i}{\Gamma_i^\ast \rho_i}$ thus $ {\Gamma_i^\ast}^2 = (\frac{{\bar{v}^\ast_i}}{{\rho^\ast_i}})^2$.
From (\ref{sys:Polar_alpha}), it holds that
$\dot{\alpha}^\ast_i = \frac{\bar{v}^\ast_i}{\rho^\ast_i}  \sin{\beta^\ast_i} = \Gamma_i^\ast \sin^2{\beta^\ast_i} = \Gamma_i^\ast$.
From Proposition \ref{pro:hat_alpha_subsystem} and Lemma \ref{th:goal_angleG_0matrix},
one can have $\hat{\alpha}^\ast = d$, $\dot{\hat{\alpha}}^\ast = \mathbf{0}_N$.
It follows that ${\Gamma_i^\ast} = \Omega$.
Since $\Gamma_i^\ast \neq 0$, the equilibrium in Case I-a only exists when $\Omega \neq 0$.
From (\ref{cal:equilibrium_pre}), one can have
$E^\ast_i  = \frac{\bar{v}^\ast_i}{ \rho^\ast_i \sin{\beta^\ast_i}} - (\frac{\bar{v}_i^\ast}{ \rho_i^\ast})^2 = \Gamma_i^\ast -{\Gamma_i^\ast}^2$.
It follows that  $\mu (\rho_i^\ast - R_i){{\rho_i^\ast}^{\sigma}} =(\Omega - \Omega^2)-(\Gamma_i^\ast -{\Gamma_i^\ast}^2)$.
Together with ${\Gamma_i^\ast} = \Omega$, one can have $\rho_i^\ast = R_i$.
Moreover, since $\Gamma_i^\ast \sin{\beta^\ast_i}>0$, i.e., $\Omega \sin{\beta^\ast_i}>0$, one can check that
$\sin{\beta^\ast_i} = 1$ for $\Omega>0$ and $\sin{\beta^\ast_i} = -1$ for $\Omega<0$.
						
To sum up, for Subcase I-a, an equilibrium (\ref{cal:equilibrium_Ia}) exists when $\Omega \neq 0$.

\chen{Subcase I-b: $\bar{v}_i ^\ast = 0, \rho_i ^\ast \neq 0, \; \forall i \in \mathcal{V}$.}
From (\ref{sys:Polar_alpha}), we get
$\ddot{\alpha}_i={\Gamma_i - \dot{\alpha}_i  + \dot{\alpha}_i  \cot{\beta_i}} -  2 \dot{\alpha}_i^2  \cot{\beta_i}$.
Since $\bar{v}_i ^\ast = 0, \rho_i ^\ast \neq 0$, one can check that $\dot{\alpha}^\ast_i=0$.
It follows that $ \ddot{\alpha}^\ast_i = \Gamma_i^\ast = 0$.
Together with the definition of $\hat{\alpha}_i$,
one can have $\dot{\hat{\alpha}}_i^\ast = 0, \ddot{\hat{\alpha}}_i^\ast = 0, \; \forall i \in \mathcal{V}$.
Thus, considering system (\ref{dynamic:hat_alpha}), the equilibrium in this case satisfies
\begin{equation}\nonumber
	 \ddot{\hat{\alpha}}^\ast 	= -\lambda_1 L(d) 	 {\hat{\alpha}}^\ast
\end{equation}
Thus, $L(d){\hat{\alpha}}^\ast = 0$. It holds that $\sum^{N}_{i=1}\hat{\alpha}_i=2\pi$ from the definition.
In view of Lemma \ref{le:eigenvalueL}, one can check that ${\hat{\alpha}}^\ast = d$.
Then we calculate $f_i$ by (\ref{sys:Polar_para1}), and get $f_i^\ast(t)=0$.
It follows $\Omega = \Gamma_i^\ast$.
Since $\Gamma_i^\ast = 0$, the equilibrium in Case I-b only exists when $\Omega = 0$.
From (\ref{sys:Polar}), we have $\ddot{\rho}_i = {  \rho_i E_i -  \dot{\alpha}_i \rho_i - \dot{\rho}_i}   + {\dot{\alpha}_i}^2  \rho_i $.
It follows $E^\ast_i=0$. Together with $\Omega = 0$, one can check that $\rho_i^\ast = R_i$.

To sum up, for Subcase I-b, an equilibrium (\ref{cal:equilibrium_Ib}) exists when $\Omega = 0$.

\chen{Subcase I-c: $\bar{v}_{i_a}^\ast \neq 0, \; \exists {i_a} \in \mathcal{V}$
	and $\bar{v}_{i_b} ^\ast = 0, \rho_{i_b} ^\ast \neq 0, \; \exists  {i_b} \in \mathcal{V}$, where $\{{i_a}\} \bigcup \{{i_b}\} = \mathcal{V}$.}
Using the similar idea with the calculation in Case I-a and I-b, one can have
\begin{small}
\begin{eqnarray}\nonumber
		\begin{cases}
			\dot{\alpha}_{i_a}^\ast = \Gamma_{i_a}^\ast \neq 0 & 	 \\
			\dot{\alpha}_{i_b}^\ast =  \Gamma_{i_b}^\ast = 0  &
		\end{cases}
\end{eqnarray}
\end{small}
It follows that
\begin{small}
\begin{eqnarray}\nonumber
		\begin{cases}
			{\alpha}_{i_a}^\ast = \Gamma_{i_a}^\ast  t + c_{i_a}& 	 \\
			{\alpha}_{i_b}^\ast = c_{i_b}   &
		\end{cases}
\end{eqnarray}
\end{small}
where $c_{i_a}, c_{i_b}$ are constants.
In this case, both ${i_a}$-agent and ${i_b}$-agent exist in the system.
It implies that there exists at least one ${i_b}$-agent (labeled as $i_b^{\prime}$) who has one or two ${i_a}$-agent as its neighbor.
One can check by (\ref{sys:Polar_para1}) that, for such an agent $i_b^{\prime}$, its $f_{i_b}^\ast$ is a function of $t$.
Thus  $\Gamma_{i_b}^\ast$ is also a function of $t$. Comparing $\Gamma_{i_b}^\ast = 0 $, we arrive at a contradiction.

To sum up, for Subcase I-c, no equilibrium exists.
\end{proof}
%

%
\keypoint{\begin{proposition}[Equilibrium Case II]\label{pro:equilibrium_caseII}
	The equilibrium of the $N$-agent system (\ref{sys:Polar}) is (\ref{cal:equilibrium_II}) for any $ \Omega$,
	if it satisfies $(\rho_i^\ast)^2 + (\bar{v}_i^\ast)^2 = 0, \; \forall i \in \mathcal{V}$.
	\begin{small}
\begin{eqnarray}\label{cal:equilibrium_II}
	(\textrm{for} \; \forall \Omega) \;  { 
		\begin{cases}
			\bar{v}_i^\ast =  0 \\
			\rho_i^\ast = 0  \\
			\dot{\alpha}_i^\ast = 0	\\
			{\hat{\alpha}}^\ast = 0			
		\end{cases} }
\end{eqnarray}
\end{small}
\end{proposition}}
%
\begin{proof}
It holds that $\rho_i^\ast=0, \bar{v}_i^\ast = 0$, since $(\rho_i^\ast)^2 + (\bar{v}_i^\ast)^2 = 0$.
Combining with  the definition of $\alpha_i$, we have $\alpha_i^\ast=0$.
Then one can check  (\ref{sys:Polar}) and (\ref{sys:Polar_alpha}) and derive that
$\dot{\alpha}_i^\ast = 0$ and thus ${\hat{\alpha}}^\ast_i = 0, \; \forall i \in \mathcal{V}$.
This completes the proof.
\end{proof}
%

%
\keypoint{\begin{proposition}[Equilibrium Case III]\label{pro:equilibrium_caseIII}
	The equilibria of the $N$-agent system (\ref{sys:Polar}) are
	(\ref{cal:equilibrium_IIIb11}) and (\ref{cal:equilibrium_IIIb2}) when $\Omega \neq 0$,
	and are  (\ref{cal:equilibrium_IIIb10}) and (\ref{cal:equilibrium_IIIb2}) when $\Omega = 0$,
	if it satisfies  $(\rho_i^\ast)^2 + (\bar{v}_i^\ast)^2 \neq 0, \;  i \in \mathcal{V}_1$
	and $(\rho_j^\ast)^2 + (\bar{v}_j^\ast)^2 = 0, \; j \in \mathcal{V}_2$, where $\mathcal V_1 \bigcup \mathcal V_2= \mathcal{V}$ and $\mathcal V_1 \bigcap \mathcal V_2=\emptyset$.
	\begin{small}
\begin{eqnarray}\label{cal:equilibrium_IIIb11}
	(\textrm{for} \; \Omega \neq 0) \;  { 
		\begin{cases}
			\keypoint{\bar{v}_j^\ast = \bar{v}_{i_b}^\ast  =  0} \\
			\keypoint{\rho_j^\ast = 0}, \; \rho_{i_b}^\ast = R_{i_b} >0\\
			\dot{\alpha}_j^\ast = \dot{\alpha}_{i_b}^\ast  = 0	\\
			\keypoint{{\hat{\alpha}}^\ast = d^{\prime} \neq d \; ({\alpha}_j^\ast = 0)}		\\
			\mbox{where}\; \revision{|\mathcal{V}_2| = 1, \; |\mathcal{V}_{1b}| = N-1	}	
		\end{cases} }
\end{eqnarray}
\begin{eqnarray}\label{cal:equilibrium_IIIb10}
	(\textrm{for} \; \Omega = 0) \;  { 
		\begin{cases}
			\bar{v}_j^\ast = \bar{v}_{i_b}^\ast  =  0 \\
			\keypoint{\rho_j^\ast = 0}, \; \rho_{i_b}^\ast = R_{i_b} >0\\
			\dot{\alpha}_j^\ast = \dot{\alpha}_{i_b}^\ast  = 0	\\
			{\hat{\alpha}}^\ast = d \; ({\alpha}_j^\ast = 0)	\\
			\mbox{where}\; |\mathcal{V}_2| = 1, \; |\mathcal{V}_{1b}| = N-1	
		\end{cases} }
\end{eqnarray}
\begin{eqnarray}\label{cal:equilibrium_IIIb2}
	(\textrm{for} \; \forall \Omega) \;  { 
		\begin{cases}
			\keypoint{\bar{v}_j^\ast = \bar{v}_{i_b}^\ast  =  0} \\
			\keypoint{\rho_j^\ast = 0}, \; \rho_{i_b}^\ast = R_{i_b} >0\\
			\dot{\alpha}_j^\ast = \dot{\alpha}_{i_b}^\ast  = 0	\\
			\keypoint{{\hat{\alpha}}^\ast = d^{\prime\prime} \neq d \; ({\alpha}_j^\ast = 0)}		\\
			\mbox{where}\; |\mathcal{V}_2| = k>1, \; |\mathcal{V}_{1b}| = N-k		
		\end{cases} }
\end{eqnarray}
\end{small}
where $i_a \in\mathcal{V}_{1a} $, $i_b \in \mathcal{V}_{1b}$,and $d^{\prime}, d^{\prime\prime}$ are constants whose value are related to the initial states.
\end{proposition}}
%
\begin{proof}
For ease of expression, we denote the agent satisfying $(\rho_j^\ast)^2 + (\bar{v}_j^\ast)^2 = 0$ by $j$-agent,
the one satisfying $(\rho_i^\ast)^2 + (\bar{v}_i^\ast)^2 \neq 0$ by $i$-agent,
the one satisfying $\bar{v}_{i_a}^\ast \neq 0$ by $i_a$-agent $i_a \in\mathcal{V}_{1a} $,
and the one satisfying $\bar{v}_{i_b} ^\ast = 0, \rho_{i_b} ^\ast \neq 0$ by $i_b$-agent $i_b \in \mathcal{V}_{1b}$.
One can check that $\mathcal{V}_1 \bigcup \mathcal{V}_2 = \mathcal{V}$,
$\mathcal{V}_1 \bigcap \mathcal{V}_2=\emptyset$,
$\mathcal{V}_{1a} \bigcup \mathcal{V}_{1b} = \mathcal{V}_1$,
and $\mathcal{V}_{1a} \bigcap \mathcal{V}_{1b} = \emptyset$,

By the same way in Proposition \ref{pro:equilibrium_caseI} and \ref{pro:equilibrium_caseII}, we have
\begin{small}
\begin{eqnarray}\nonumber	
			j\textrm{-agent}: \; (\rho_j^\ast)^2 + (\bar{v}_j^\ast)^2 = 0  {\Rightarrow}   %
				\begin{cases} {\alpha}_j^\ast = 0\\  \dot{\alpha}_j^\ast= 0 \end{cases}\\\nonumber
			{{i_a\textrm{-agent}: \; \bar{v}_{i_a}^\ast \neq 0}   {\Rightarrow}  }
				\begin{cases} {\alpha}_{i_a}^\ast = \Gamma_{i_a}^\ast t + c_{i_a} \\ \dot{\alpha}_{i_a}^\ast=\Gamma_{i_a}^\ast \neq 0 \end{cases}\\\nonumber
			i_b\textrm{-agent}: \; \bar{v}_{i_b} ^\ast = 0, \; \rho_{i_b} ^\ast \neq 0   {\Rightarrow}
				\begin{cases} {\alpha}_{i_b}^\ast = c_{i_b} \\ \dot{\alpha}_{i_b}^\ast= 0 \end{cases} 			
\end{eqnarray}
\end{small}
where $c_{i_a}, c_{i_b}$ are constants.

It's worth to point that, this case requires $\mathcal{V}_1  \neq \emptyset , \mathcal{V}_2 \neq \emptyset$.
Together with $\mathcal{V}_{1a} \bigcup \mathcal{V}_{1b} = \mathcal{V}_1$, in this case, we need to consider further two subcases.
\chen{Subcase III-a: $\mathcal{V}_{1a}  \neq \emptyset$}.
\chen{Subcase III-b: $\mathcal{V}_{1a}  = \emptyset, \mathcal{V}_{1b}  \neq \emptyset$}.
Due to a space limit, we omit the details,
which are similar to those involved in the proofs of  Proposition \ref{pro:equilibrium_caseI} and Proposition \ref{pro:equilibrium_caseII},
and present directly the results as follows.

For Subcase III-a, no equilibrium exists.
For Subcase III-b, when $|\mathcal{V}_2| = 1, \; |\mathcal{V}_{1b}| = N-1$,
equilibrium (\ref{cal:equilibrium_IIIb11}) exists for $\Omega \neq 0$ and (\ref{cal:equilibrium_IIIb10}) exists for $\Omega = 0$;
while when $ |\mathcal{V}_2| = k>1, \; |\mathcal{V}_{1b}| = N-k$, equilibrium (\ref{cal:equilibrium_IIIb2}) for $\forall \Omega$.
This completes the proof.
\end{proof}
%

\revision{Note that the equilibrium (\ref{cal:equilibrium_Ia}) or (\ref{cal:equilibrium_Ib}) corresponds to the desired formation structure. }
Finally, we summarize the above discussion on the  equilibria of the $N$-agent system (\ref{sys:Polar}) in the following proposition.
%
\keypoint{\begin{proposition}\label{cal:equilibrium_all}
	The $N$-agent system (\ref{sys:Polar}) has
	equilibria (\ref{cal:equilibrium_Ia}), (\ref{cal:equilibrium_IIIb11}), (\ref{cal:equilibrium_II}), (\ref{cal:equilibrium_IIIb2}) when $\Omega \neq 0$,
	and has equilibria (\ref{cal:equilibrium_Ib}), (\ref{cal:equilibrium_IIIb10}), (\ref{cal:equilibrium_II}), (\ref{cal:equilibrium_IIIb2}) when $\Omega = 0$.
\end{proposition}}
%
%

\subsection{Analysis of convergence}


Following the idea used to design the controller,
we first focus on the case when only the first part, \revision{the} converging part, of the controller works.

We emphasize that the converging part concerns the relationship between each agent and the \revision{static/}moving target, and no interaction between agents is included.
Thus, to analyze the system's convergence under the converging part of the controller, we just need to focus on each agent and the target,
and therefore it can be regarded as a single-agent \cut{case.} 
It's obvious that for this case, the system has equilibria (\ref{cal:equilibrium_Ia}), (\ref{cal:equilibrium_Ib}), and (\ref{cal:equilibrium_II}).
%
\keypoint{\begin{theorem}
Suppose that $N=1$. The equilibrium (\ref{cal:equilibrium_Ia}) corresponding to $\Omega\neq0$ and the equilibrium (\ref{cal:equilibrium_Ib}) corresponding to $\Omega=0$ are locally stable equilibria. When $\sigma=0,$ the equilibrium (\ref{cal:equilibrium_II}) is unstable.
\end{theorem}}
%
\begin{proof}
For the single-agent case, the system (\ref{sys:Polar}) becomes
\begin{eqnarray}\nonumber
	\begin{cases}
		\dot{\rho} = \bar{v} \cos{\beta}  \\
		\dot{\bar{v}} =  { \rho (   E \cos{\beta} + \Omega \sin{\beta}) -  \bar{v} } \\
		\dot{\beta} = { 1 - \frac{\rho}{\bar{v}} (   E \sin{\beta} - \Omega  \cos{\beta} ) } - \frac{\bar{v}}{\rho}\sin{\beta}
 	\end{cases}
\end{eqnarray}
where $E= -\mu (\rho - R){\rho^{\sigma}}  - {\Omega}(\Omega - 1)$.

For $\Omega\neq0,$ its Jacobian matrix denoted as $A$ can be derived as
\begin{small}
\begin{equation*}\nonumber
\begin{split}
	A|_{(\ref{cal:equilibrium_Ia})}
= \left[\begin{array}{ccc}
			0 & 0 & - \Omega R \\
			\Omega \sin{\beta} & -1 & R \Omega (\Omega-1) \sin{\beta}  \\
			\frac{2 \Omega -1}{R} + \mu  \frac{R^{\sigma}}{\Omega} & \frac{\sin{\beta}}{R}(\frac{1}{\Omega}-2) & -1
		\end{array}\right].
\end{split}
\end{equation*}
\end{small}
Its characteristic polynomial is given by
\begin{small}
\begin{equation*}\nonumber
\begin{split}
	& \det (\lambda I - A|_{(\ref{cal:equilibrium_Ia})})  \\\nonumber
		&= \lambda^3 + 2\lambda^2 +[(2\Omega - 1)^2+1 + \mu  R^{\sigma+1}]\lambda + \mu  R^{\sigma+1}.
\end{split}
\end{equation*}
\end{small}
It is easy to see  that the coefficients of the polynomial  $a_0\triangleq1$, $a_1\triangleq2$, $a_2\triangleq(2\Omega - 1)^2+1 + \mu  R^{\sigma+1}$, $a_3\triangleq\mu  R^{\sigma+1}$ are all positive. In addition, one can verify that $a_1 a_2  - a_0 a_3>0$. Therefore, $\det (\lambda I - A|_{(\ref{cal:equilibrium_Ia})})$ is stable and the equilibrium (\ref{cal:equilibrium_Ia}) is locally stable.

When $\Omega=0$, the Jacobian matrix at the equilibrium (\ref{cal:equilibrium_Ib}) can be calculated as
\begin{small}
\begin{eqnarray}\nonumber
	A|_{(\ref{cal:equilibrium_Ib})}= \left[\begin{array}{ccc}
			0 & \cos{\beta} & 0 \\
			- \mu  R^{\sigma+1} \cos{\beta}   & -1 & 0 \\
			\ast & \ast & -1
		\end{array}\right].
\end{eqnarray}
\end{small}
The elements denoted by $\ast$ is irrelevant for the calculation of the characteristic polynomial. The characteristic polynomial of $A|_{(\ref{cal:equilibrium_Ib})}$  is given by
\begin{small}
\begin{eqnarray}\nonumber
	 \det (\lambda I - A|_{(\ref{cal:equilibrium_Ib})})  = (\lambda+1)[\lambda^2+\lambda+\mu  R^{\sigma+1} (\cos{\beta})^2].
\end{eqnarray}
\end{small}
Since at the equilibrium, $\beta=0$ and hence $\mu  R^{\sigma+1} (\cos{\beta})^2>0$, the eigenvalues of the characteristic polynomial are in the open left half plane and the equilibrium (\ref{cal:equilibrium_Ib}) is locally stable.

When $\sigma=0,$  the equilibrium (\ref{cal:equilibrium_II}) of system (\ref{sys:Polar}) corresponds to the following equilibrium $({\bar{x}}^\ast,{\bar{y}}^\ast,{\bar{v}}^{x^\ast},{\bar{v}}^{y^\ast})=(0,0,0,0)$
of system (\ref{sys:XY_1}), which is the single-agent case of the $N$-agent system (\ref{dynamic:double_i}) under control laws (\ref{control:part1}).
\begin{small}
\begin{eqnarray}\label{sys:XY_1}
	\begin{cases}
		\dot{\bar{x}} = \bar{v}^x  \\
		\dot{\bar{y}} = \bar{v}^y  \\
		\dot{\bar{v}}^x =  \revision{\bar{x} E} - {\Omega} \bar{y} - \bar{v}^y -  \bar{v}^x \\  
		\dot{\bar{v}}^y =  \revision{\bar{y} E} + {\Omega} \bar{x} + \bar{v}^x - \bar{v}^y
	\end{cases}
\end{eqnarray}
\end{small}
where $E= -\mu (\rho - R){\rho^{\sigma}}  - {\Omega}(\Omega - 1)$.
Its Jacobian matrix at this equilibrium is given by
\begin{small}
\begin{equation*}\nonumber
\begin{split}
	A|_{(\ref{cal:equilibrium_II})}
= \left[\begin{array}{cccc}
			0 & 0 & 1 & 0 \\
			0 & 0 & 0 & 1 \\
			\mu R - \Omega(\Omega-1) & - \Omega& -1 & -1 \\
			\Omega & \mu R - \Omega(\Omega-1) & 1 & -1
		\end{array}\right],
\end{split}
\end{equation*}
\end{small}
whose characteristic polynomial can be calculated as
\begin{small}
\begin{eqnarray}\nonumber
	 \det (\lambda I - A|_{(\ref{cal:equilibrium_II})})
		=& [ \lambda(\lambda+1) - \mu R + \Omega(\Omega-1) ]^2 +(\lambda+\Omega)^2 \\\nonumber
		=& \lambda^4 + 2 \lambda^3 + 2(\Omega^2-\Omega-\mu R +1)\lambda^2 \\\nonumber
		& + 2(\Omega^2-\mu R)\lambda + [\Omega(\Omega-1) -\mu R ]^2 + \Omega^2
\end{eqnarray}
\end{small}
The Routh array of this polynomial is given as
%
\begin{table}[H] \nonumber
	\begin{tabular}{m{0.1cm}|m{3cm}|m{1.8cm}|m{1.8cm}}
		\hline		
		$\lambda^4$ & $1$ & $2(\Omega^2-\Omega-\mu R +1)$ & $[\Omega(\Omega-1) -\mu R ]^2 + \Omega^2$ \\
		\hline
		$\lambda^3$ & $2$ & $2(\Omega^2-\mu R)$ & $0$ \\
		\hline
		$\lambda^2$ & $\Omega^2-2\Omega-\mu R +2$ &  $[\Omega(\Omega-1) -\mu R ]^2 + \Omega^2$ &  \\
		\hline
		$\lambda^1$ & $\frac{-4\mu R }{\Omega^2-2\Omega-\mu R +2}$ & $0$ &  \\
	 	\hline
	  	$\lambda^0$ & $[\Omega(\Omega-1) -\mu R ]^2 + \Omega^2  $ &  &  \\
		\hline
	\end{tabular}
\end{table}\label{table1}
%
%

If $\Omega^2-2\Omega-\mu R +2>0$ (resp. $<0$), then $\frac{-4\mu R}{\Omega^2-2\Omega-\mu R +2}<0$  (resp. $>0$).
In both cases, the number of sign changes in the first column of the array is two, so the characteristic polynomial has two eigenvalues with positive real parts and therefore is unstable.  If $\Omega^2-2\Omega-\mu R +2 = 0$, one can also show that characteristic polynomial also has two eigenvalues with positive real parts. One concludes that the  equilibrium (\ref{cal:equilibrium_II}) is unstable.
\end{proof}
%

For the case when the layout part is included and then the integrated controller is considered, the analysis becomes more challenging and we do not have a complete result yet.
Nevertheless, we will demonstrate the effectiveness  of our controller by simulations in the next section.
\par~

%
\section{Simulation results}\label{se:simulations}
%
%
In the simulations, we consider a system consisting of six agents.
The target starts from \revision{or stays still at} the point $(0,0)$ in the plane without loss of generality.
The initial states of the six agents \revision{are} generated randomly.

We present \revision{three} typical scenarios.
\revision{For Example $1$, the desired general formation is a circle rotating clockwise around a static target,
where the agents' rotating velocity $\Omega=-0.2$,
the distance from agents to the target is $1$,
and the angular distances between neighbors are set to be equal.}
\revision{For Example $2$}, the desired general formation is a form of two concentric circles rotating counterclockwise around \revision{a moving} target,
where the agents' rotating velocity $\Omega=1$,
the distance from \revision{the} agent to the target is $0.6$ for agent $1,3,5$ and $1.5$ for agent $2,4,6$,
and the angular distances between neighbors are set to be equal.
\revision{For Example $3$}, the desired general formation is a \revision{right} triangle remaining static relative to \revision{a moving} target,
where the agents' rotating velocity \cut{$\Omega=0$.}
%
\cut{For} ease of comparison, for each case, we use the same parameters of the controllers as $\lambda_1=\lambda_2=1, \mu=1, \sigma=-1$;
\revision{and the same trajectories of the moving target in Example $2$ and $3$.}

\begin{figure*}[thpb]   
\begin{center}
           \subfigure[\revision{Example 1 (static target): $\Omega=-0.2$}]{\label{fig:0a}
           \includegraphics[width=0.36\linewidth]{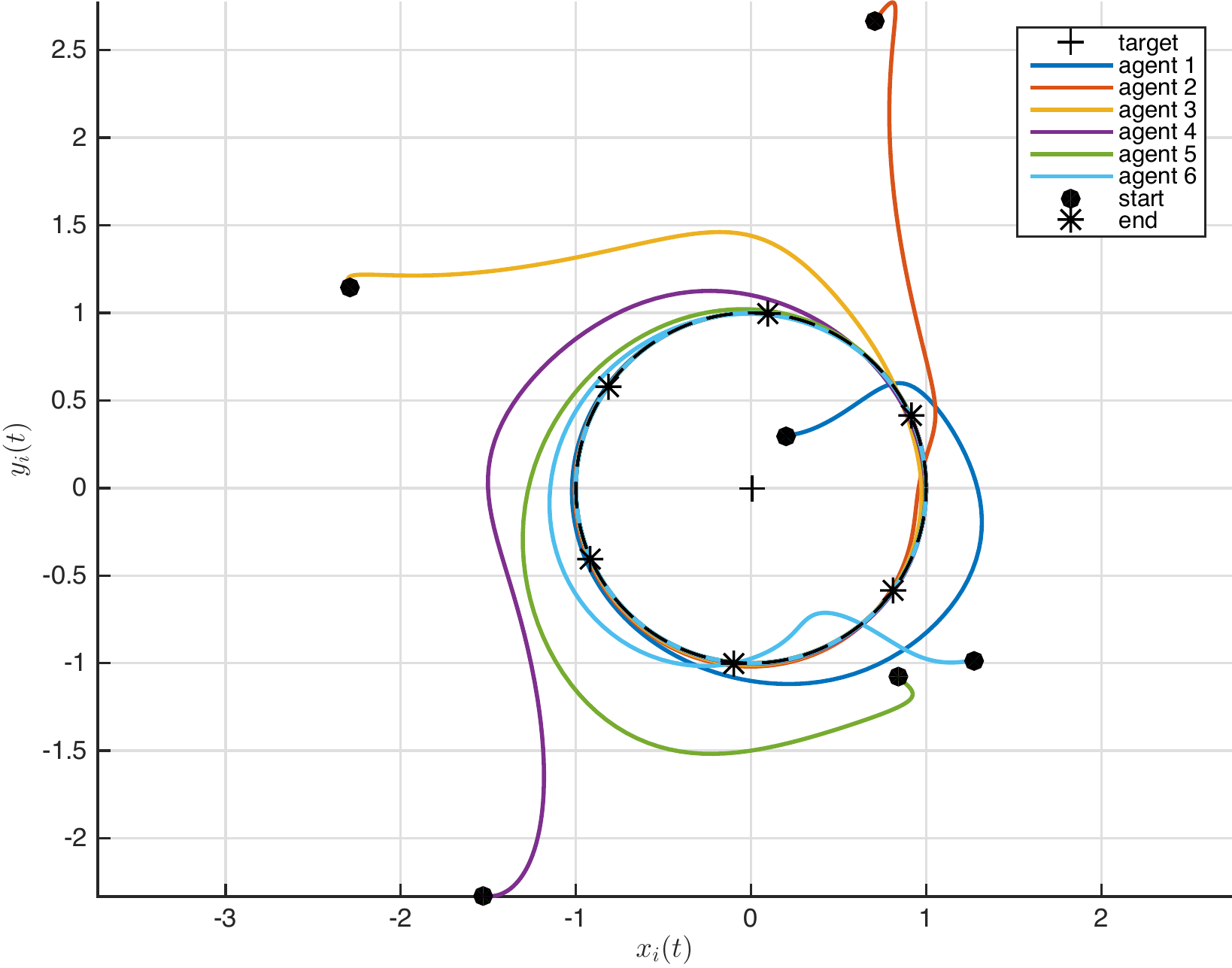}}
           \subfigure[\revision{Example 1 (static target): $\Omega=-0.2$}]{\label{fig:0b}
           \includegraphics[width=0.36\linewidth]{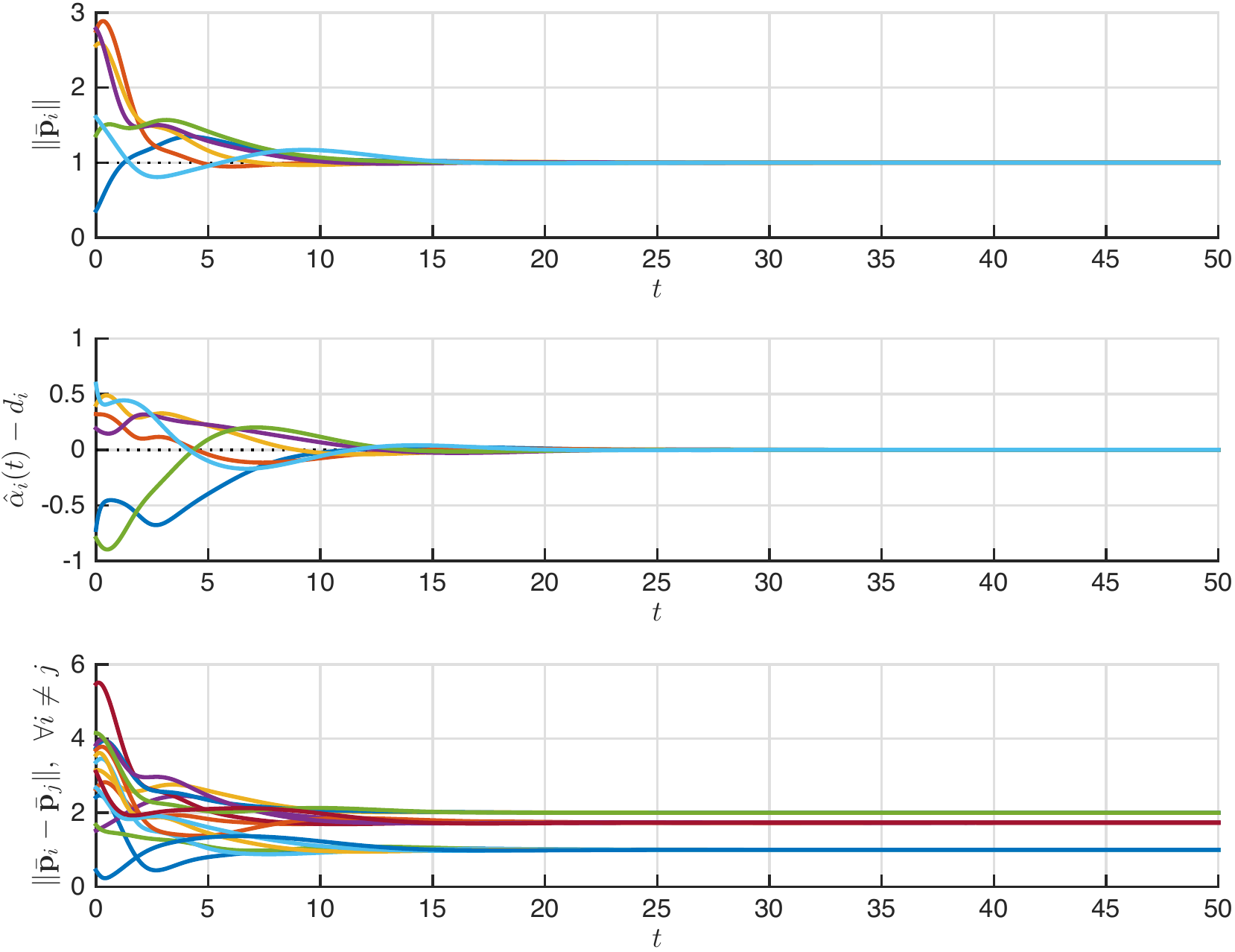}}
           \subfigure[{Example \revision{2 (moving target)}: $\Omega=1$}]{\label{fig:1a}
           \includegraphics[width=0.36\linewidth]{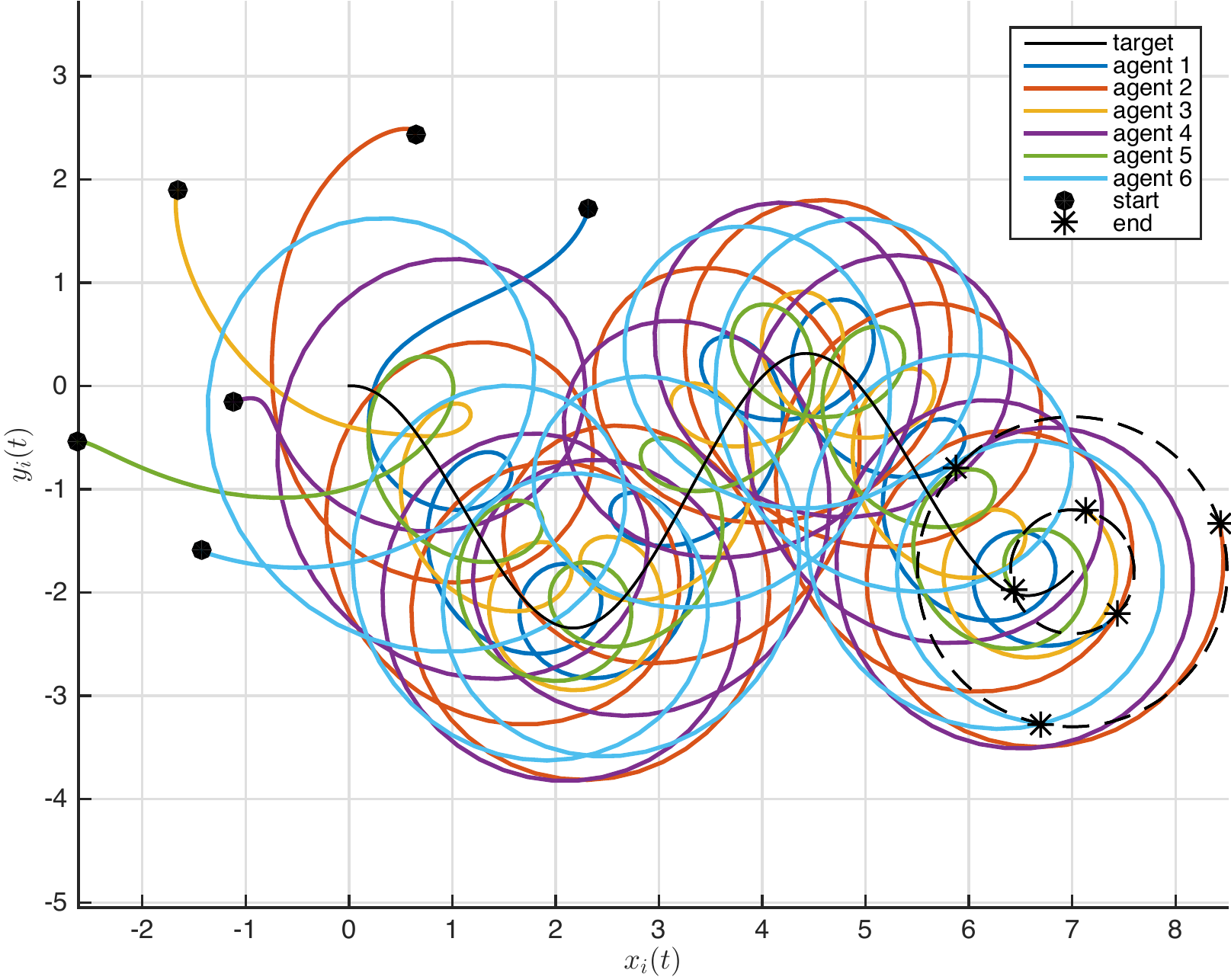}}
           \subfigure[{Example \revision{2 (moving target)}: $\Omega=1$}]{\label{fig:1b}
           \includegraphics[width=0.36\linewidth]{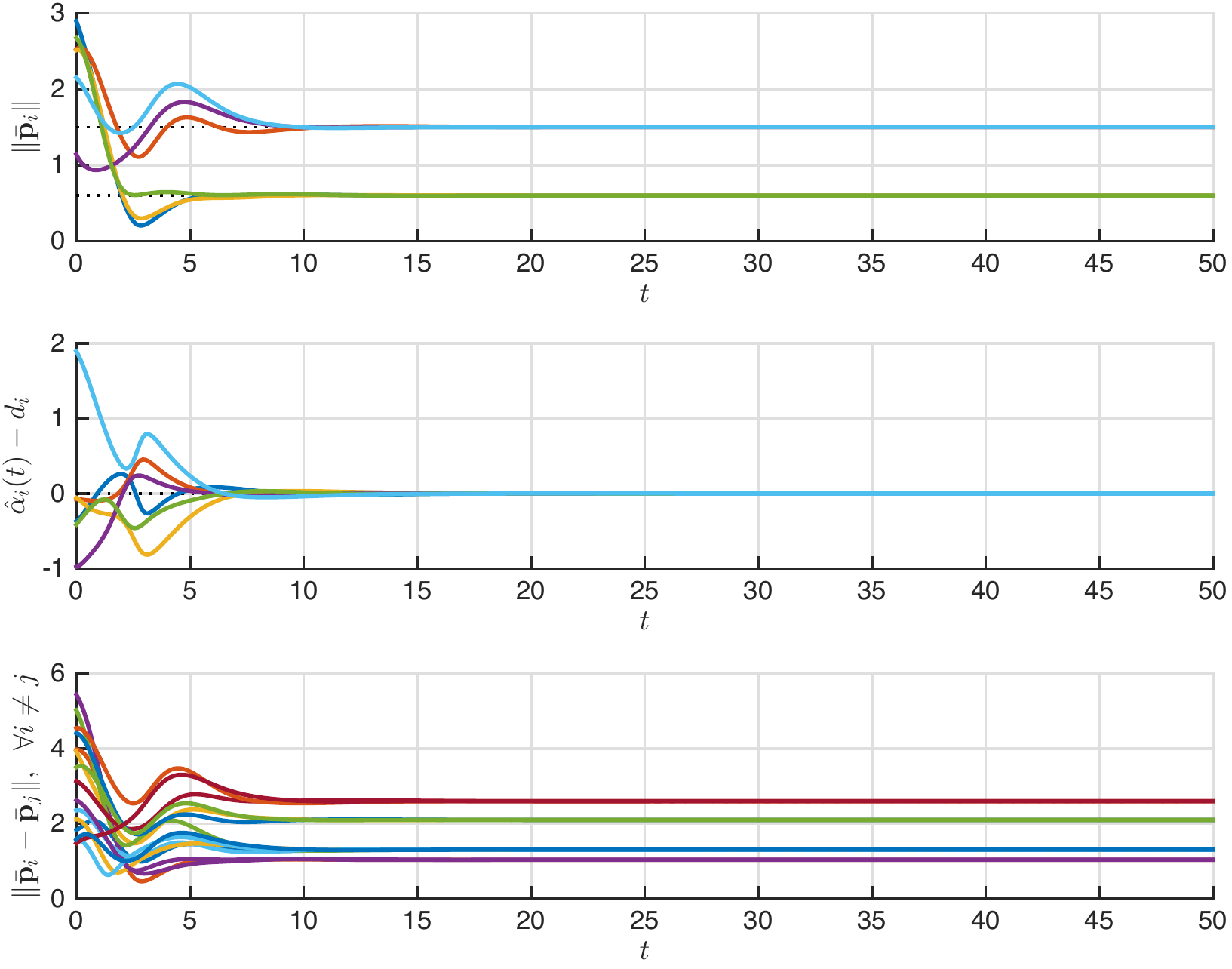}}
           \subfigure[{Example \revision{3 (moving target)}: $\Omega=0$}]{\label{fig:2a}
           \includegraphics[width=0.36\linewidth]{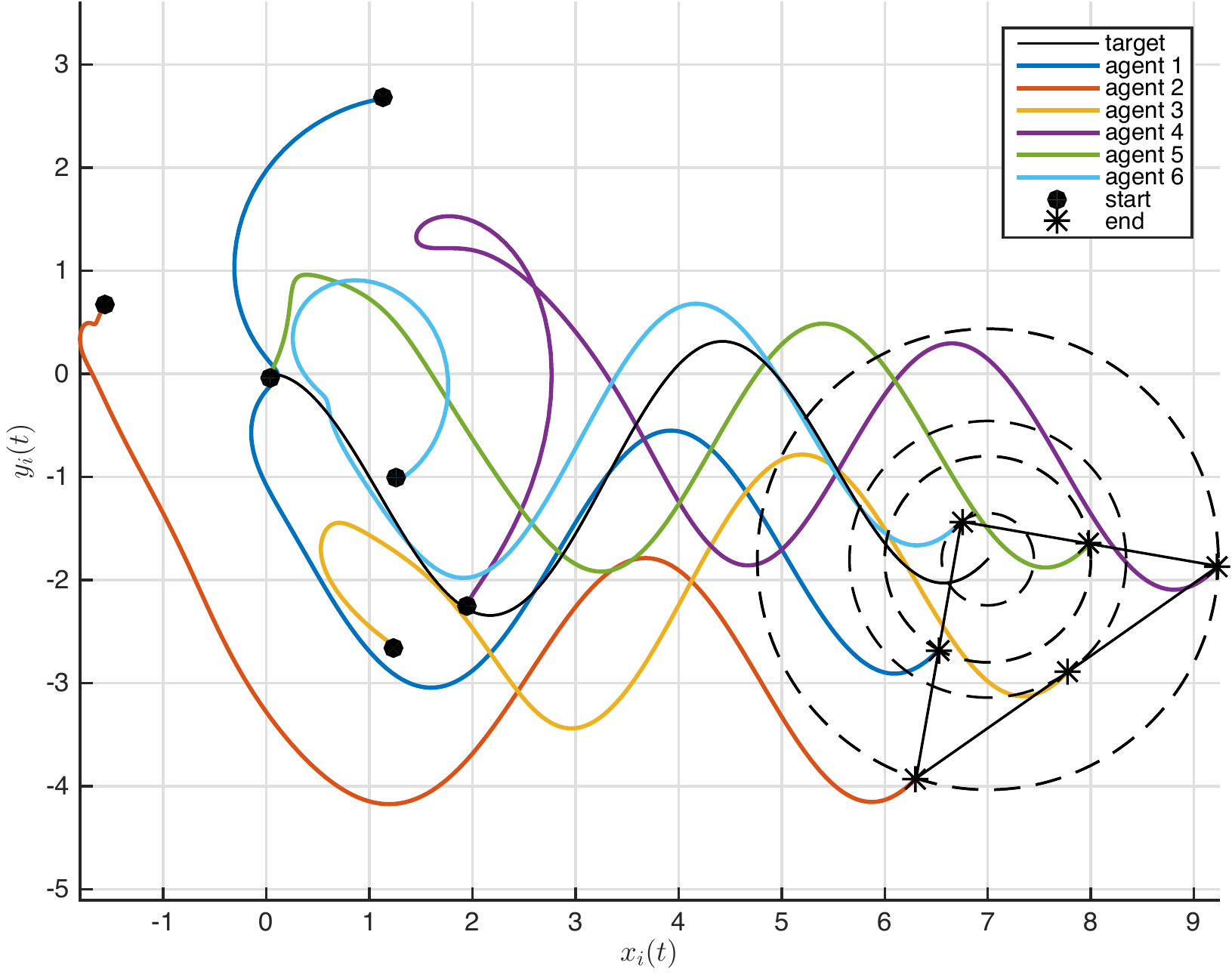}}
           \subfigure[{Example \revision{3 (moving target)}: $\Omega=0$}]{\label{fig:2b}
           \includegraphics[width=0.36\linewidth]{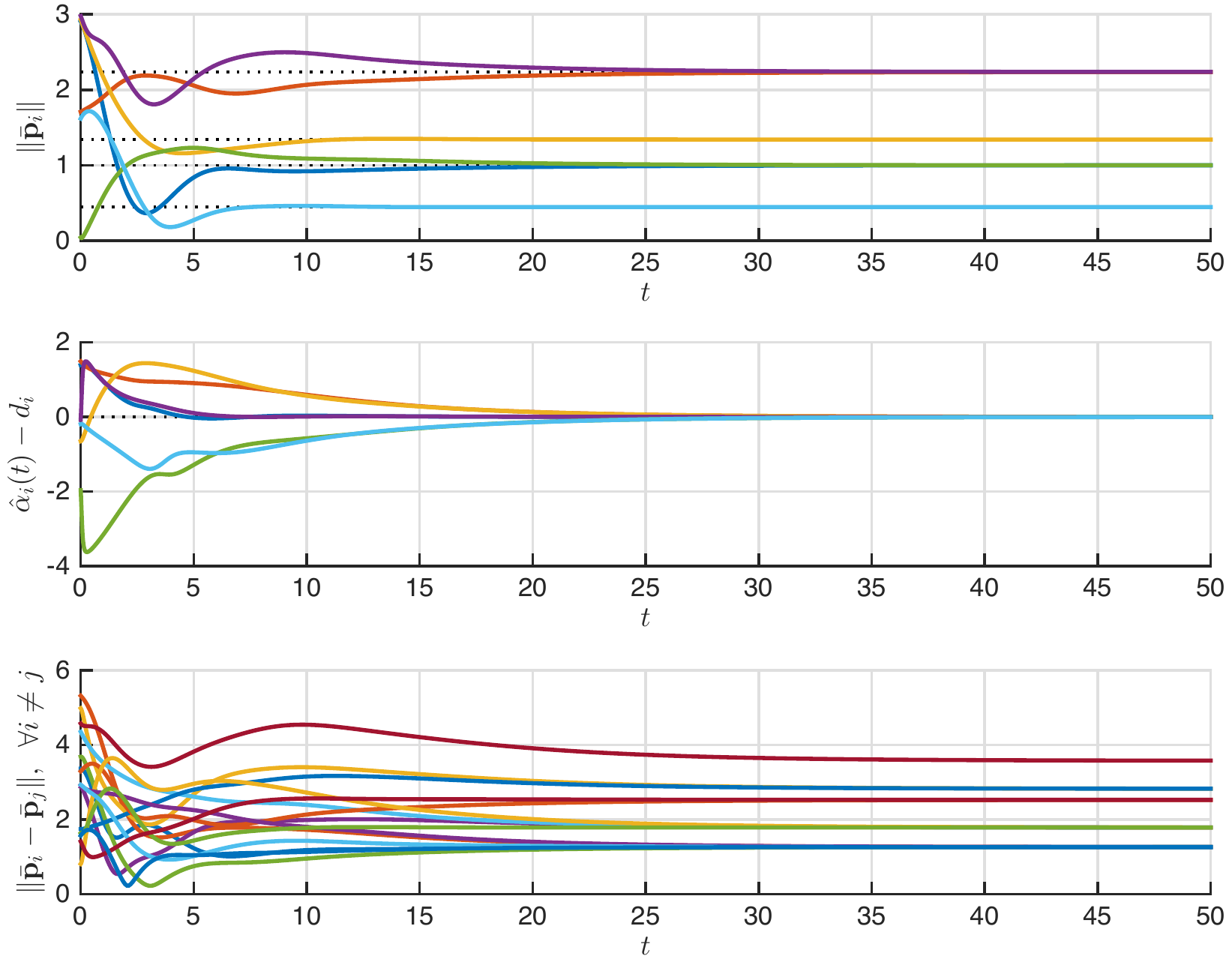}}
   \caption{Simulation results of the proposed controllers solving the general formation problem when
   	$N=6$, and $\lambda_1=\lambda_2=1, \mu=1, \sigma=-1$.
	\revision{(a)(b) Example 1 with a static target: $\Omega=-0.2$.}
	\revision{(c)(d)} Example \revision{2 with a moving target}: $\Omega=1$.
	\revision{(e)(f)} Example \revision{3 with a moving target}: $\Omega=0$.
	\revision{(a)(c)(e)} Trajectories of six agents in the plane;
	\revision{(b)(d)(f)} Distances from each agent to the target,
	differences between current distances and the desired distances between all pairs of neighbors,
	and distances between all pairs of agents.
   }\label{fig:example}
\end{center} 
\end{figure*}

We run the simulations and show the results in \cut{Fig. \ref{fig:example}.}
%
%
\cut{The} simulation results clearly indicate \revision{that} our proposed controllers solve the general formation problem, while no collision ever occurs.

\par~

%
\section{Conclusions}\label{se:conclusions}
%
In this paper, we have studied the general formation problem for a group of mobile agents in a plane.
The problem includes two \cut{sub-objectives: target circling and distribution adjustment.}
%
%
\cut{Using} the limit-cycle based design idea, we have designed a distributed local controller combined \revision{two} parts to solve the general control \cut{problem.}
%
\cut{The} theoretical analysis has provided to show the convergence of the system under \revision{the} converging part of the proposed controller,
while numerical simulations have been performed to demonstrate the performance of the whole controller.
Currently, we are working on the theoretical analysis when the layout part is \revision{included in the whole controller.}

\par~
\par~

%
%
%
%
\bibliographystyle{IEEEtran}        

\end{document}